\documentclass[11pt]{article}
\usepackage{cite}
\usepackage{amssymb}
\usepackage{amsfonts}
\usepackage{mathrsfs}
\usepackage{graphicx}
\usepackage{amsmath}
\usepackage{color}
\usepackage{amsthm}
\usepackage{epstopdf}
\usepackage{booktabs}
\usepackage{pifont,bm}
\usepackage{tikz}
\usepackage{enumerate}
\usepackage{threeparttable}
\usepackage{diagbox}

\theoremstyle{definition}

\makeatletter
\def\@biblabel#1{[#1]}
\makeatother

\makeatletter \@addtoreset{equation}{section}

\begin{document}

\begin{titlepage}
\title{\bf{Long-time asymptotics for the reverse space-time nonlocal Hirota equation with decaying initial value problem: Without solitons
\footnote{
Corresponding authors.\protect\\
\hspace*{3ex} E-mail addresses: ychen@sei.ecnu.edu.cn (Y. Chen)}
}}
\author{Wei-Qi Peng$^{a}$, Yong Chen$^{a,b,*}$\\
\small \emph{$^{a}$School of Mathematical Sciences, Shanghai Key Laboratory of PMMP} \\
\small \emph{East China Normal University, Shanghai, 200241, China} \\
\small \emph{$^{b}$College of Mathematics and Systems Science, Shandong University }\\
\small \emph{of Science and Technology, Qingdao, 266590, China} \\
\date{}}
\thispagestyle{empty}
\end{titlepage}
\maketitle

\vspace{-0.5cm}
\begin{center}
\rule{15cm}{1pt}\vspace{0.3cm}

\parbox{15cm}{\small
{\bf Abstract}\\
\hspace{0.5cm}  In this work, we mainly consider the Cauchy problem for the reverse space-time nonlocal Hirota equation with the initial data rapidly decaying in the solitonless sector. Start from the Lax pair,  we first construct the basis Riemann-Hilbert problem for the reverse space-time nonlocal Hirota equation. Furthermore, using the approach of Deift-Zhou  nonlinear steepest descent,  the explicit long-time asymptotics  for the reverse space-time nonlocal Hirota is derived. For the reverse space-time nonlocal Hirota equation, since the symmetries of its scattering matrix are different with the local Hirota equation, the $\vartheta(\lambda_{i})(i=0, 1)$ would like to be imaginary, which results in the $\delta_{\lambda_{i}}^{0}$ contains an increasing  $t^{\frac{\pm Im\vartheta(\lambda_{i})}{2}}$, and then the asymptotic behavior for nonlocal Hirota equation becomes differently.
}

\vspace{0.5cm}
\parbox{15cm}{\small{

\vspace{0.3cm} \emph{Key words:} Riemann-Hilbert problem; Reverse space-time nonlocal Hirota equation; Long-time asymptotics.\\

 } }
\end{center}
\vspace{0.3cm} \rule{15cm}{1pt} \vspace{0.2cm}

\section{Introduction}
In recent years, more and more scholars pay attention to nonlocal integrable equations in the area of integrable systems. The nonlocal nonlinear Schr\"{o}dinger (NLS) equation
\begin{align}\label{0.1}
iq_{t}+q_{xx}+2q^{2}q^{\ast}(-x, t)=0
\end{align}
was first introduced by Ablowitz and
Musslimani, and they derived its soliton solutions through the method of inverse scattering transform(IST)\cite{Tian-CPAM1,Tian-CPAM2}. The nonlocal NLS equation\eqref{0.1} contains the $PT$ symmetric potential which is invariant under the transformation $x\rightarrow -x$ and complex conjugation. The IST for the nonlocal NLS equation with nonzero boundary conditions at infinity was studied in Ref.\cite{JMP-S2}. It is worth mentioning that the long-time asymptotics for the integrable nonlocal NLS equation with decaying  boundary conditions has been presented in Ref.\cite{JMP-S}. Moreover, other nonlocal integrable equations were also investigated including nonlocal Davey-Stewartson equations, nonlocal modified KdV equation, nonlocal sine-Gordon equation, nonlocal derivative NLS equation, etc. \cite{zhang-33,zhang-34,zhang-35,zhang-36, zhang-38,zhang-39}.

Recently,  the reverse space-time nonlocal Hirota equation
\begin{align}\label{0.2}
iq_{t}+\alpha\left[q_{xx}-2q^{2}q(-x, -t)\right]+i\beta\left[q_{xxx}-6qq(-x,-t)q_{x}\right]=0,
\end{align}
was introduced by Cen, Correa and Fring in 2019 \cite{Tian-CPAM8} as a $PT$ symmetric reduction ($r=q(-x, -t)$) of the following system:
\begin{align}\label{2}
iq_{t}+\alpha\left[q_{xx}-2q^{2}r\right]+i\beta\left[q_{xxx}-6qrq_{x}\right]=0,\notag\\
ir_{t}-\alpha\left[r_{xx}-2qr^{2}\right]+i\beta\left[r_{xxx}-6qrr_{x}\right]=0,
\end{align}
where $\alpha, \beta\in \mathbb{R}$. These two equations are transformed into each other by means of a $PT$-symmetry transformation. This reduction leads to \eqref{0.2} and is consistent with the $PT$-symmetry condition: if $q(x, t)$ is a solution of \eqref{1}, then $q(-x, -t)$ is a solution as well. The other types of integrable nonlocal Hirota equation are also introduced in Ref. \cite{Tian-CPAM8} including
the reverse space nonlocal Hirota equation, reverse time nonlocal Hirota equation  and
the conjugate reverse space-time nonlocal Hirota equation. The multi-soliton solutions of these nonlocal Hirota equations have been generated by DarbouxCrum transformations and Hirota bilinear method\cite{Tian-CPAM8}.  The conjugate reverse space-time nonlocal Hirota equation with nonzero boundary conditions was investigated via Riemann-Hilbert(RH) method\cite{Peng-PD}.  Soliton solutions  of the conjugate reverse space nonlocal Hirota equation were obtained by the IST method and Darboux transformation method in Ref. \cite{Tian-CPAM} and Ref. \cite{ND-guo}, respectively. Using the Darboux transformation,  some types of exact solutions of the reverse space-time nonlocal Hirota equation were found in Ref. \cite{CPB-xia}.

In 1967, the IST was used to solve the KdV equation with Lax pairs by Gardner et al. for the first time\cite{Wangdengshan-na2}. Since then, it played an increasingly important role in finding the exact solutions for integrable systems. Later on, a modern version of IST method, named RH
formulation, was developed by Zakharov et al.\cite{Wangdengshan-na20}, and then the exact
solutions and long-time asymptotics of various integrable equations were investigated using RH formulation \cite{Wangdengshan-na21,Wangdengshan-na22,Wangdengshan-na23,Wangdengshan-na24, Wangdengshan-na25}. It is worth mentioning that Pelinovsky and Shimabukuro proved the existence of global solutions to the derivative
NLS equation on the line from the perspective of inverse scattering transform based on the
representation of a RH problem, which is a milestone in the development of IST\cite{Pelinovsky}.   And what's more, in 1993, Deift-Zhou put forward the nonlinear steepest descent method for the first time to solve the oscillatory RH problem and derive the long-time asymptotics of solutions for the modified KdV equation\cite{DeiftZ}. After that,  this method has been employed
to discuss the asymptotic analysis in a variety of  integrable models\cite{Peng16-13,Peng16-14,Peng16-15,Peng16-16,Peng16-17,Peng16-18,Peng16-19,Peng16-20}. In 2019, Dmitry Shepelsky et al. applied this method to study the long-time behavior of solutions to the initial boundary value problem of nonlocal NLS equations\cite{JMP-S}.  Recently, they have extended the Deift-Zhou method to study the long-time asymptotic behavior of nonlocal integrable NLS solutions with nonzero boundary conditions and step-like initial data, respectively\cite{Peng16-21,Peng16-22}. Besides, the Deift-Zhou nonlinear steepest-descent method was used to analyze the long-time asymptotics for the solution of the nonlocal mKdV equation\cite{Peng16-23}.

As we know, long-time asymptotics for the reverse space-time nonlocal
Hirota equation \eqref{0.2} has not been reported. In this paper,
we are committed to the Cauchy problem for the so-called defocusing reverse space-time nonlocal Hirota equation
\begin{align}\label{1}
&iq_{t}+\alpha\left[q_{xx}-2q^{2}q(-x, -t)\right]+i\beta\left[q_{xxx}-6qq(-x,-t)q_{x}\right]=0, \quad x\in \mathbb{R}, \quad t>0, \notag\\
&q(x,0)=q_{0}(x),
\end{align}
with the initial data $q_{0}(x)$  belonging  to the Schwartz space and rapidly decaying to $0$ as $\mid x\mid\rightarrow \infty$. For the nonlocal Hirota equation, except to the symmetries of its Lax pair are different with the local Hirota equation, another major difference is that $\vartheta(\lambda_{i})(i=0, 1)$ are imaginary in our case, which leads to the $\delta_{\lambda_{i}}^{0}$ contains an increasing  $t^{\frac{\pm Im\vartheta(\lambda_{i})}{2}}$, and then the asymptotic behavior for nonlocal Hirota equation will behave differently. The pivotal result of this paper is generalized in what follows:

\noindent \textbf{Theorem 1.1.} \emph{
Let $q(x, t)$ be the solution of the Cauchy problem of the reverse space-time nonlocal Hirota equation\eqref{1} with $q_{0}(x)$  lying in the Schwartz space. As $t\rightarrow\infty$, for $\alpha^{2}-3\beta \frac{t}{x}>0$, the leading asymptotics of the solution $q(x, t)$ is
\begin{gather}
q(x, t)=\frac{\sqrt{\pi}t^{-\frac{1}{2}-Im\vartheta(\lambda_{0})}
e^{4i\lambda_{0}^{2}t(4\beta\lambda_{0}+\alpha)+2\chi_{0}(\lambda_{0})+\frac{\pi\vartheta(\lambda_{0}) }{2}+\frac{\pi i}{4}+iRe\vartheta(\lambda_{0})\ln t+i\vartheta(\lambda_{0})\ln(32\lambda_{0}^{2}(\alpha+6\beta\lambda_{0}))}}{\sqrt{\alpha+6\beta\lambda_{0}}r_{1}(\lambda_{0})\Gamma(i\vartheta(\lambda_{0}))}\notag\\
+\frac{\sqrt{\pi}t^{-\frac{1}{2}+Im\vartheta(\lambda_{1})}
e^{4i\lambda_{1}^{2}t(4\beta\lambda_{1}+\alpha)+2\chi_{1}(\lambda_{1})-\frac{\pi}{2}\vartheta(\lambda_{1})+\frac{\pi i}{4}
-iRe\vartheta(\lambda_{1})\ln t-i\vartheta(\lambda_{1})\ln(32\lambda_{1}^{2}(\alpha+6\beta\lambda_{1})) }}{\sqrt{\alpha+6\beta\lambda_{1}}r_{1}(\lambda_{1})\Gamma(-i\vartheta(\lambda_{1}))}\notag\\
+O(t^{-\frac{1}{2}-\mbox{max}\{|\mbox{Im}\vartheta(\lambda_{0})|,|\mbox{Im}\vartheta(\lambda_{1})|\}}),
\end{gather}
with
\begin{gather}
\lambda_{0}=\frac{-\alpha-\sqrt{\alpha^{2}-3\beta\frac{x}{t}}}{6\beta},\ \lambda_{1}=\frac{-\alpha+\sqrt{\alpha^{2}-3\beta \frac{x}{t}}}{6\beta},\notag\\
\chi_{0}(\lambda)=\frac{1}{2\pi i}\int_{\lambda_{0}}^{\lambda_{1}}\ln\left(\frac{1-r_{1}(s)r_{2}(s)}
{1-r_{1}(\lambda_{0})r_{2}(\lambda_{0})}\right)\frac{\mathrm{d}s}{s-\lambda},\notag\\
\chi_{1}(\lambda)=\frac{1}{2\pi i}\int_{\lambda_{0}}^{\lambda_{1}}\ln\left(\frac{1-r_{1}(s)r_{2}(s)}
{1-r_{1}(\lambda_{1})r_{2}(\lambda_{1})}\right)\frac{\mathrm{d}s}{s-\lambda},\notag\\
\vartheta(\lambda_{0})=-\frac{1}{2\pi}\ln(1-r_{1}(\lambda_{0})r_{2}(\lambda_{0})),\notag\\
\vartheta(\lambda_{1})=-\frac{1}{2\pi}\ln(1-r_{1}(\lambda_{1})r_{2}(\lambda_{1})),
\end{gather}
where $\Gamma$ is Gamma function, and $r_{1}(\lambda), r_{2}(\lambda)$ are the reflection coefficients.}

\textbf{Organization of the paper:} In section 2,  we perform the direct scattering theory to generate the associated RH problem, further the phase analysis is discussed in detail. In section 3, the nonlinear steepest descent method is used to analyse the long-time asymptotics of the solution for the reverse space-time nonlocal Hirota equation.

\section{Inverse scattering transform and the Riemann-Hilbert problem}
At the very start, we should carry out the direct   scattering analysis  to construct the basis RH problem for the nonlocal Hirota equation\eqref{1}, which is a compatibility condition of the following Lax pair\cite{Tian-CPAM8}
\begin{gather}
\phi_{x}=X\phi,\qquad X\equiv\left(\begin{array}{cc}
    -i\lambda  &  q(x,t) \\
    q(-x, -t)  &  i\lambda \\
\end{array}\right),\notag\\
\phi_{t}=T\phi,\qquad T\equiv \left(\begin{array}{cc}
    Q  &  B \\
    C  &  -Q \\
\end{array}\right),\notag\\
Q=-i\alpha qq(-x, -t)-2i\alpha\lambda^{2}+\beta\left[q(-x, -t)q_{x}+qq_{x}(-x, -t)-4i\lambda^{3}-2i\lambda qq(-x, -t)\right],\notag\\
B=i\alpha q_{x}+2\alpha\lambda q+\beta\left[2q^{2}q(-x, -t)-q_{xx}+2i\lambda q_{x}+4\lambda^{2}q\right],\notag\\
C=i\alpha q_{x}(-x, -t)+2\alpha\lambda q(-x, -t)\notag\\
+\beta\left[2qq(-x, -t)^{2}-q_{xx}(-x, -t)+2i\lambda q_{x}(-x, -t)+4\lambda^{2}q(-x, -t)\right],\label{4}
\end{gather}
where $\lambda$ means the spectral parameter, $\phi=\phi\left(x, t; \lambda\right)$ is the eigenfunction.

As $x\rightarrow \pm \infty$, due to the initial data rapidly decaying, the Lax pair \eqref{4} turns into
\begin{align}\label{5}
\phi^{\infty}_{x}=X_{0}\phi^{\infty}=-i\lambda\sigma_{3}\phi^{\infty}, \qquad \phi^{\infty}_{t}=T_{0}\phi^{\infty}=(2\alpha\lambda+4\beta\lambda^{2})X_{0}\phi^{\infty},
\end{align}
where $\sigma_{3}$ represents one of the following Pauli matrices
\begin{align}
\sigma_{1}=\left(\begin{array}{cc}
    0  &  1\\
    1 &  0\\
\end{array}\right), \quad \sigma_{2}=\left(\begin{array}{cc}
    0  &  -i\\
    i &  0\\
\end{array}\right), \quad \sigma_{3}=\left(\begin{array}{cc}
    1  &  0\\
    0 &  -1\\
\end{array}\right).
\end{align}
We can easily find that the Eq.\eqref{5} arrives at the following fundamental matrix solution
\begin{align}\label{6}
\phi^{\infty}(x, t; \lambda)=e^{-i\theta(x, t; \lambda)\sigma_{3}},\qquad \theta(x, t; \lambda)=\lambda(x+\left[2\alpha\lambda+4\beta\lambda^{2}\right]t).
\end{align}
Taking $D_{+}$, $D_{-}$ and $\Sigma$ on $\lambda$-plane as $D_{\pm}= \left\{\lambda\in \mathbb{C} | \mbox{Im}\lambda\gtrless 0\right\}, \Sigma=\mathbb{R}$, the Jost solutions $\phi_{\pm}(x, t; \lambda)$ are
\begin{align}\label{7}
\phi_{\pm}(x, t; \lambda)=e^{-i\theta(x, t; \lambda)\sigma_{3}}+o(1), \quad \lambda\in \Sigma, \quad \mbox{as} \quad x\rightarrow \pm \infty.
\end{align}
Through the variable transformation
\begin{align}\label{8}
\mu_{\pm}(x, t; \lambda)=\phi_{\pm}(x, t; \lambda)e^{i\theta(x, t; \lambda)\sigma_{3}},
\end{align}
the spectral problem \eqref{4} can be solved as the following Jost solutions $\mu_{\pm}$, given by
\begin{align}\label{9}
\left\{
\begin{array}{lr}
\mu_{-}=I+\int_{-\infty}^{x}\exp\left[-i\lambda\sigma_{3}\left(x-y\right)\right](X-X_{0})\mu_{-}
\exp\left[i\lambda\sigma_{3}\left(x-y\right)\right]\mathrm{d}y, \\
\mu_{+}=I-\int_{x}^{+\infty}\exp\left[-i\lambda\sigma_{3}\left(x-y\right)\right](X-X_{0})\mu_{+}
\exp\left[i\lambda\sigma_{3}\left(x-y\right)\right]\mathrm{d}y.
  \end{array}
\right.
\end{align}
\noindent \textbf{Proposition 2.1.} \emph{
Suppose $q\in L^{1} (\mathbb{R}^{\pm})$, then $\mu_{\pm}(x,t,\lambda)$ given in Eq.\eqref{9} satisfy the following properties:}

\emph{$\bullet$ $\mu_{-1}(x, t, \lambda)$ and $\mu_{+2}(x, t, \lambda)$ is analytical in $D_{+}$ and
continuous in $D_{+}\cup \Sigma$,}

\emph{$\bullet$ $\mu_{+1}(x, t, \lambda)$ and $\mu_{-2}(x, t, \lambda)$  is analytical in $D_{-}$ and
continuous in $D_{-}\cup \Sigma$,}

\emph{$\bullet$ $\mu_{\pm}(x,t,\lambda)\rightarrow I$ \mbox{as}  $\lambda\rightarrow \infty$,}

\emph{$\bullet$ $\det \mu_{\pm}(x, t, \lambda)=1, \quad x, t\in \mathbb{R}, \quad \lambda\in \Sigma$.}

Since the Jost solutions $\phi_{\pm}(x, t, \lambda)$ are the simultaneous solutions of Lax pair \eqref{4}, they can meet the following linear relation contacted by a scattering matrix $S(\lambda)=(s_{i j} (\lambda))_{2\times 2}$, given by
\begin{align}\label{10}
\phi_{-}(x, t, \lambda)=\phi_{+}(x, t, \lambda)S(\lambda), \quad \lambda\in \Sigma,
\end{align}
of which the scattering coefficients can be written into what follows
\begin{align}\label{11}
s_{11}(\lambda)=Wr(\phi_{-1},\phi_{+2}), \quad s_{12}(\lambda)=Wr(\phi_{-2},\phi_{+2}),\notag\\
s_{21}(\lambda)=Wr(\phi_{+1},\phi_{-1}), \quad s_{22}(\lambda)=Wr(\phi_{+1},\phi_{-2}).
\end{align}

\noindent \textbf{Proposition 2.2.} \emph{ Suppose $q\in L^{1} (\mathbb{R}^{\pm})$, then
the scattering  matrix $S(\lambda)$ has the following properties:
}

\emph{$\bullet$ $\det S(\lambda)=1$ for $\lambda\in \Sigma$,}

\emph{$\bullet$ $s_{11}(\lambda)$ is analytical in $D_{+}$ and
continuous in $D_{+}\cup \Sigma$,}

\emph{$\bullet$ $s_{22}(\lambda)$ is analytical in $D_{-}$ and
continuous in $D_{-}\cup \Sigma$,}

\emph{$\bullet$ $S(x,t,\lambda)\rightarrow I$ \mbox{as}  $\lambda\rightarrow \infty$.}

Furthermore, we need to study the symmetries of the Jost solutions $\phi(x, t, \lambda)$ and scattering matrix $S(\lambda)$ for the nonlocal Hirota equation. The detail reduction conditions for $X(x, t, \lambda)$ and $T(x, t, \lambda)$ in the Lax pair \eqref{4} on $\lambda$-plane are as follows:
\begin{align}\label{12}
X(x, t, \lambda)=-\sigma_{2}X(-x, -t, \lambda)\sigma_{2},\qquad T(x, t, \lambda)=-\sigma_{2}T(-x, -t, \lambda)\sigma_{2},
\end{align}
which results in the Jost solutions $\Psi(x, t, \lambda)$, and scattering matrix $S(\lambda)$ has the following
reduction conditions on $\lambda$-plane:
\begin{align}\label{13}
\phi_{\pm}(x, t, \lambda)=\sigma_{2}\phi_{\mp}(-x, -t, \lambda)\sigma_{2},\qquad S(\lambda)=\sigma_{2}S^{-1}(\lambda)\sigma_{2},
\end{align}
which means $s_{12}(\lambda)=s_{21}(\lambda)$, and $s_{11}(\lambda), s_{22}(\lambda)$ are not directly related. This is different from the  case of local Hirota equation.

According to the analyticity of Jost solutions $\mu_{\pm}(x, t, \lambda)$ in Proposition 2.1, we can define the following sectionally meromorphic matrices
\begin{align}\label{14}
M_{+}(x, t, \lambda)=(\frac{\mu_{-1}}{s_{11}},\mu_{+2}),\qquad M_{-}(x, t, \lambda)=(\mu_{+1},\frac{\mu_{-2}}{s_{22}}),
\end{align}
where $\pm$ denote analyticity in $D_{+}$ and $D_{-}$, respectively. Then, a matrix RH problem is generated:

\noindent \textbf{Riemann-Hilbert Problem}  \emph{
$M(x, t, \lambda)$ solves the following RH problem:
\begin{align}\label{15}
\left\{
\begin{array}{lr}
M(x, t, \lambda)\ \mbox{is analytic in} \ \mathbb{C }\setminus \Sigma,\\
M_{+}(x, t, \lambda)=M_{-}(x, t, \lambda)J(x, t, \lambda), \qquad \lambda\in \Sigma,\\
M(x, t, \lambda)\rightarrow I,\qquad \lambda\rightarrow \infty,\\
  \end{array}
\right.
\end{align}
of which the jump matrix $J(x, t, \lambda)$ is
\begin{align}\label{16}
J(x, t, \lambda)=\left(\begin{array}{cc}
    1-r_{1}(\lambda)r_{2}(\lambda)  &  -r_{2}(\lambda)e^{-2i\theta(x, t, \lambda)}\\
  r_{1}(\lambda)e^{2i\theta(x, t, \lambda)} &  1\\
\end{array}\right),
\end{align}
where $r_{1}(\lambda)=\frac{s_{21}(\lambda)}{s_{11}(\lambda)}, r_{2}(\lambda)=\frac{s_{12}(\lambda)}{s_{22}(\lambda)}$.
}

Let
\begin{align}\label{17}
M(x, t, \lambda)=I+\frac{1}{\lambda}M_{1}(x, t; \lambda)+O(\frac{1}{\lambda^{2}}),\qquad \lambda\rightarrow \infty,
\end{align}
then the potential $q(x, t)$ of the nonlocal Hirota equation \eqref{1} is given by
\begin{align}\label{18}
q(x, t)=2i\left[M_{1}\right]_{12}(x, t, \lambda)=2i\lim_{\lambda\rightarrow\infty}\lambda \left[M\right]_{12}(x, t, \lambda).
\end{align}

\section{The long-time behavior for the nonlocal Hirota equation}
In this section, we primarily devote to discuss the long-time behavior for the nonlocal Hirota equation\eqref{1}.
Let's start with phase analysis, in terms of the works of Deift and Zhou \cite{DeiftZ}, we take $\frac{df}{\mathrm{d}\lambda}=0$, and then the stationary points of the function $f$ are $\lambda_{0}=\frac{-\alpha-\sqrt{\alpha^{2}-3\beta \xi}}{6\beta}, \lambda_{1}=\frac{-\alpha+\sqrt{\alpha^{2}-3\beta \xi}}{6\beta}$ for $\alpha^{2}-3\beta \xi>0$, there we have defined $f=\lambda(\xi+2\alpha\lambda+4\beta\lambda^{2}), \xi=\frac{x}{t}$, and the signature distribution for $Re (i f)$ is shown in Fig. 1. The steepest decent contours are
\begin{align}\label{19}
&L:\{\lambda=\lambda_{1}+\lambda_{1}\rho e^{\frac{3\pi i}{4}}:-\infty<\rho\leq\sqrt{2}\}
\cup \{\lambda=\lambda_{0}-\lambda_{0}\rho e^{\frac{\pi i}{4}}:-\infty<\rho\leq\sqrt{2}\},\notag\\
&L^{\ast}:\{\lambda=\lambda_{1}+\lambda_{1}\rho e^{-\frac{3\pi i}{4}}:-\infty<\rho\leq\sqrt{2}\}
\cup \{\lambda=\lambda_{0}-\lambda_{0}\rho e^{-\frac{\pi i}{4}}:-\infty<\rho\leq\sqrt{2}\}.
\end{align}

\centerline{\begin{tikzpicture}[scale=1.0]
\draw[-][thick](-3,0)--(-2,0);
\draw[-][thick](-2.0,0)--(-1.0,0)node[below]{$\lambda_{0}$};
\draw[-][thick](-1,0)--(0,0);
\draw[-][thick](0,0)--(1,0)node[below]{$\lambda_{1}$};
\draw[-][thick](1,0)--(2,0);
\draw[fill] (1,0) circle [radius=0.035];
\draw[fill] (-1,0) circle [radius=0.035];
\draw[-][thick](2.0,0)--(3.0,0);
\draw [-,thick, cyan] (1,0) to [out=90,in=-120] (1.9,3);
\draw [-,thick, cyan] (1,0) to [out=-90,in=120] (1.9,-3);
\draw [-,thick, cyan] (-1,0) to [out=90,in=-60] (-1.9,3);
\draw [-,thick, cyan] (-1,0) to [out=-90,in=60] (-1.9,-3);
\draw[fill] (0,2) node{$Re(if)>0$};
\draw[fill] (0,-2) node{$Re(if)<0$};
\draw[fill] (-2,1) node{$Re(if)<0$};
\draw[fill] (2,1) node{$Re(if)<0$};
\draw[fill] (2,-1) node{$Re(if)>0$};
\draw[fill] (-2,-1) node{$Re(if)>0$};
\end{tikzpicture}}
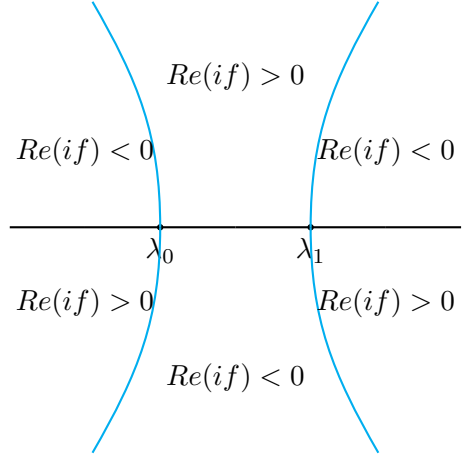
\centerline{\noindent {\small \textbf{Figure 1.} (Color online) The signature table for $Re( if)$ in the complex $\lambda$-plane.}}

\subsection{Factorization of the jump matrix and contour deformation}
We decompose the jump matrix $J(x, t, \lambda)$ into following two cases:
\begin{align}\label{20}
J(x, t, \lambda)=\left\{
\begin{array}{lr}
\left(\begin{array}{cc}
     1 &  -r_{2}(\lambda)e^{-2ift} \\
     0  &  1\\
\end{array}\right)\left(\begin{array}{cc}
     1 &  0\\
     r_{1}(\lambda)e^{2ift}  &  1\\
\end{array}\right),\\
\left(\begin{array}{cc}
     1 &  0\\
     \frac{r_{1}(\lambda)e^{2ift}}{1-r_{1}(\lambda)r_{2}(\lambda)}  &  1\\
\end{array}\right)\left(\begin{array}{cc}
     1-r_{1}(\lambda)r_{2}(\lambda) &  0 \\
     0  &  \frac{1}{1-r_{1}(\lambda)r_{2}(\lambda)}\\
\end{array}\right)
\left(\begin{array}{cc}
    1 &  -\frac{r_{2}(\lambda)e^{-2ift}}{1-r_{1}(\lambda)r_{2}(\lambda)}\\
    0  &  1\\
\end{array}\right).
  \end{array}
\right.
\end{align}
Then, we define a  RH problem for the function $\delta(\lambda)$
\begin{align}\label{21}
\left\{
\begin{array}{lr}
\delta_{+}( \lambda)=(1-r_{1}(\lambda)r_{2}(\lambda))\delta_{-}( \lambda), \qquad \lambda\in(\lambda_{1},\lambda_{2}),\\
\delta(\lambda)\rightarrow 1,\qquad \lambda\rightarrow \infty,
  \end{array}
\right.
\end{align}
which can be solved by the Plemelj formula as
\begin{align}\label{22}
\delta(\lambda)=\exp\left\{\frac{1}{2\pi i}\int_{\lambda_{0}}^{\lambda_{1}}\frac{\ln(1-r_{1}(s)r_{2}(s))}{s-\lambda}\mathrm{d}s\right\}
=\left(\frac{\lambda-\lambda_{1}}{\lambda-\lambda_{0}}\right)^{i\vartheta(\lambda_{0})}e^{\chi_{0}(\lambda)}
=\left(\frac{\lambda-\lambda_{1}}{\lambda-\lambda_{0}}\right)^{i\vartheta(\lambda_{1})}e^{\chi_{1}(\lambda)},
\end{align}
where
\begin{align}\label{23}
\chi_{0}(\lambda)=\frac{1}{2\pi i}\int_{\lambda_{0}}^{\lambda_{1}}\ln\left(\frac{1-r_{1}(s)r_{2}(s)}
{1-r_{1}(\lambda_{0})r_{2}(\lambda_{0})}\right)\frac{\mathrm{d}s}{s-\lambda},\notag\\
\chi_{1}(\lambda)=\frac{1}{2\pi i}\int_{\lambda_{0}}^{\lambda_{1}}\ln\left(\frac{1-r_{1}(s)r_{2}(s)}
{1-r_{1}(\lambda_{1})r_{2}(\lambda_{1})}\right)\frac{\mathrm{d}s}{s-\lambda},
\end{align}
\begin{align}\label{24}
\vartheta(\lambda_{0})=-\frac{1}{2\pi}\ln(1-r_{1}(\lambda_{0})r_{2}(\lambda_{0})),\notag\\
\vartheta(\lambda_{1})=-\frac{1}{2\pi}\ln(1-r_{1}(\lambda_{1})r_{2}(\lambda_{1})),
\end{align}
so that
\begin{align}\label{25}
\mbox{Im}\vartheta(\lambda_{i})=-\frac{1}{2\pi}\int_{-\infty}^{\lambda_{i}}\mathrm{d}\ \mbox{arg} (1-r_{1}(s)r_{2}(s)),\ i=0, 1.
\end{align}
Assuming that $\int_{-\infty}^{\lambda_{i}}\mathrm{d}\ \mbox{arg} (1-r_{1}(s)r_{2}(s))\in (-\pi, \pi)$, we have
\begin{align}\label{26}
|\mbox{Im}\vartheta(\lambda)|<\frac{1}{2}, \qquad \lambda\in \mathbb{R},
\end{align}
then we get that $\ln(1-r_{1}(\lambda)r_{2}(\lambda))$ is single-valued, and the singularity of $\delta(\lambda, \xi)$ at $\lambda=\lambda_{0}$ and $\lambda=\lambda_{1}$ is square integrable.

Let
\begin{align}\label{27}
M^{(1)}(x, t;\lambda)=M(x, t;\lambda)\delta^{-\sigma_{3}}(\lambda),
\end{align}
then $M^{(1)}$ solves the  RH problem on the jump contour $\mathbb{R}$ shown in Fig. 2,
\begin{align}\label{28}
\left\{
\begin{array}{lr}
M^{(1)}_{+}(x, t; \lambda)=M^{(1)}_{-}(x, t; \lambda)J^{(1)}(x, t; \lambda), \qquad \lambda\in\mathbb{R}=\Sigma^{(1)},\\
M^{(1)}(x, t; \lambda)\rightarrow I,\qquad \lambda\rightarrow \infty,
  \end{array}
\right.
\end{align}
where
\begin{align}\label{29}
J^{(1)}=
\left(\begin{array}{cc}
     1 &  0\\
     \gamma_{1}(\lambda)e^{2ift}\delta_{-}^{-2}  &  1\\
\end{array}\right)
\left(\begin{array}{cc}
    1 &  -\gamma_{2}(\lambda)e^{-2ift}\delta_{+}^{2}\\
    0  &  1\\
\end{array}\right),
\end{align}
the  functions $\gamma_{1}(\lambda),\gamma_{2}(\lambda)$ are defined as
\begin{align}\label{30}
\gamma_{1}(\lambda)=
\begin{cases}
\frac{r_{1}(\lambda)}{1-r_{1}(\lambda)r_{2}(\lambda)},\quad \lambda_{0}<\lambda<\lambda_{1},\\
-r_{1}(\lambda),\quad \lambda<\lambda_{0}\cup \lambda>\lambda_{1},\\
\end{cases}\quad  \gamma_{2}(\lambda)=\begin{cases}
\frac{r_{2}(\lambda)}{1-r_{1}(\lambda)r_{2}(\lambda)},\quad \lambda_{0}<\lambda<\lambda_{1},\\
-r_{2}(\lambda),\quad \lambda<\lambda_{0}\cup \lambda>\lambda_{1}.\\
\end{cases}
\end{align}
\\

\centerline{\begin{tikzpicture}[scale=1.5]
\draw[-][thick](-3,0)--(-2,0);
\draw[<-][thick](-2.0,0)--(-1.0,0)node[below]{$\lambda_{0}$};
\draw[->][thick](-1,0)--(0,0);
\draw[-][thick](0,0)--(1,0)node[below]{$\lambda_{1}$};
\draw[-][thick](1,0)--(2,0);
\draw[fill] (1,0) circle [radius=0.035];
\draw[fill] (-1,0) circle [radius=0.035];
\draw[<-][thick](2.0,0)--(3.0,0)node[right]{$\mathbb{R}$};
\end{tikzpicture}}
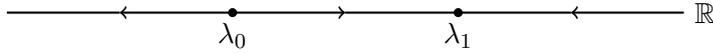
\centerline{\noindent {\small \textbf{Figure 2.} (Color online) The jump contour $\mathbb{R}=\Sigma^{(1)}$.}}

Performing the decomposition $J^{(1)}=(b_{-})^{-1}b_{+}$, where
\begin{align}\label{31}
b_{-}=\left(\begin{array}{cc}
     1 &  0\\
     -\gamma_{1}(\lambda)e^{2ift}\delta_{-}^{-2}  &  1\\
\end{array}\right),\quad
b_{+}=\left(\begin{array}{cc}
    1 &  -\gamma_{2}(\lambda)e^{-2ift}\delta_{+}^{2}\\
    0  &  1\\
\end{array}\right),
\end{align}
and taking
\begin{align}\label{32}
M^{(2)}=\left\{
\begin{array}{lr}
M^{(1)}(\lambda), \qquad \qquad \quad \lambda\in \Omega_{1}\cup \Omega_{2},\\
M^{(1)}(\lambda)(b_{-})^{-1}, \qquad \lambda\in \Omega_{3}\cup \Omega_{4}\cup\Omega_{5},\\
M^{(1)}(\lambda)(b_{+})^{-1}, \qquad \lambda\in \Omega_{6}\cup \Omega_{7}\cup\Omega_{8},\\
  \end{array}
\right.
\end{align}
we can deform the contour $\Sigma^{(1)}$ into the contour $\Sigma^{(2)}=L\cup L^{\ast}$ as displayed in Fig. 3 and derive the following RH problem on the contour $\Sigma^{(2)}=L\cup L^{\ast}\cup \mathbb{R}$
\begin{align}\label{33}
\left\{
\begin{array}{lr}
M^{(2)}_{+}(x, t; \lambda)=M^{(2)}_{-}(x, t; \lambda)J^{(2)}(x, t; \lambda), \qquad \lambda\in\Sigma^{(2)},\\
M^{(2)}(x, t; \lambda)\rightarrow I,\qquad \lambda\rightarrow \infty,
  \end{array}
\right.
\end{align}
where  the jump matrix is
\begin{align}\label{34}
J^{(2)}=\left\{
\begin{array}{lr}
I, \qquad \qquad \lambda\in\mathbb{R},\\
b_{+},  \qquad\quad\ \lambda\in L,\\
(b_{-})^{-1}, \qquad \lambda\in L^{\ast}.
  \end{array}
\right.
\end{align}

\centerline{\begin{tikzpicture}[scale=2.0]
\draw[-][thick](-3,0)--(-2,0);
\draw[<-][thick](-2.0,0)--(-1.0,0)node[below]{$\lambda_{0}$};
\draw[->][thick](-1,0)--(0,0);
\draw[->][thick](0,1)--(0.5,0.5);
\draw[-][thick](0.5,0.5)--(1,0);
\draw[->][thick](0,-1)--(0.5,-0.5);
\draw[-][thick](0.5,-0.5)--(1,0);
\draw[->][thick](2,-1)--(1.5,-0.5);
\draw[-][thick](1.5,-0.5)--(1,0);
\draw[->][thick](2,1)--(1.5,0.5);
\draw[-][thick](1.5,0.5)--(1,0);
\draw[-][thick](0,1)--(-0.5,0.5);
\draw[->][thick](-1,0)--(-0.5,0.5);
\draw[-][thick](0,-1)--(-0.5,-0.5);
\draw[->][thick](-1,0)--(-0.5,-0.5);
\draw[->][thick](-1,0)--(-1.5,0.5);
\draw[-][thick](-1.5,0.5)--(-2,1);
\draw[->][thick](-1,0)--(-1.5,-0.5);
\draw[-][thick](-1.5,-0.5)--(-2,-1);
\draw[-][thick](0,0)--(1,0)node[below]{$\lambda_{1}$};
\draw[-][thick](1,0)--(2,0);
\draw[fill] (1,0) circle [radius=0.035];
\draw[fill] (-1,0) circle [radius=0.035];
\draw[<-][thick](2.0,0)--(3.0,0);
\draw[fill] (0,1.2) node{$\Omega_{1}$};
\draw[fill] (0,-1.2) node{$\Omega_{2}$};
\draw[fill] (0,0.5) node{$\Omega_{7}$};
\draw[fill] (0,-0.5) node{$\Omega_{4}$};
\draw[fill] (2,-0.5) node{$\Omega_{6}$};
\draw[fill] (2,0.5) node{$\Omega_{3}$};
\draw[fill] (-2,-0.5) node{$\Omega_{8}$};
\draw[fill] (-2,0.5) node{$\Omega_{5}$};
\end{tikzpicture}}
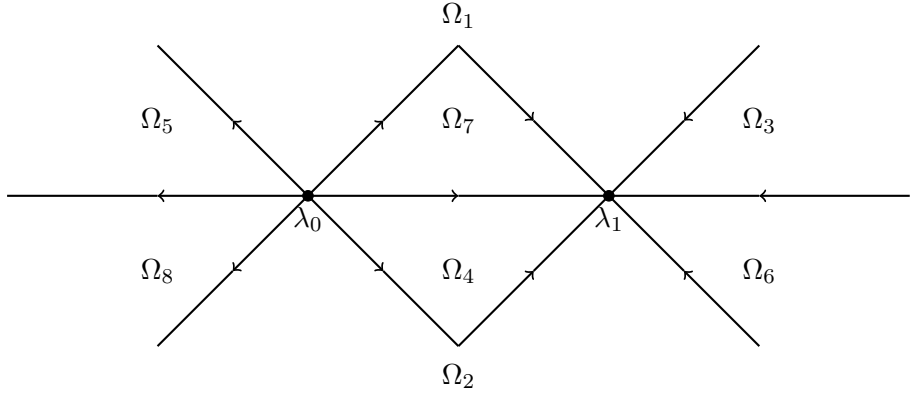
\centerline{\noindent { \textbf{Figure 3.} (Color online) The jump contour $\Sigma^{(2)}$  and domains $\Omega_{j}(j=1,\cdots, 8)$.}}

Considering the jump matrix $J^{(2)}$ decaying exponentially to identity away from the stationary phase point $\lambda_{0}, \lambda_{1}$ as $t\rightarrow\infty$,
we need take $D_{\lambda_{0}}^{\epsilon}$ and $D_{\lambda_{1}}^{\epsilon}$  be a disk of radius $\epsilon$ centered at $\lambda_{0}$ and $\lambda_{1}$, with $\epsilon$ sufficiently small. Thus, we can change the contour $\Sigma^{(2)}$ into the contours $\Sigma^{(app)}$
and $\Sigma^{(err)}$ (see Fig. 4).

Define
\begin{align}\label{35}
M^{(app)}=M_{\lambda_{0}}^{(app)}M_{\lambda_{1}}^{(app)}=\left\{
\begin{array}{lr}
I, \qquad \qquad \qquad \qquad \mbox{outside}  \quad D_{\lambda_{0}}^{\epsilon}\cup D_{\lambda_{1}}^{\epsilon},\\
\mbox{parametrix of} \ M^{(2)}, \quad \mbox{inside} \quad D_{\lambda_{0}}^{\epsilon}\cup D_{\lambda_{1}}^{\epsilon},
  \end{array}
\right.
\end{align}
which means $M^{(app)}$  has the same jump conditions as $M^{(2)}$ inside $D_{\lambda_{0}}^{\epsilon}\cup D_{\lambda_{1}}^{\epsilon}$. $M_{\lambda_{0}}^{(app)}$ should possess a jump
$J_{\lambda_{0}}^{(app)}$ across the circle $D_{\lambda_{0}}^{\epsilon}$, $M_{\lambda_{1}}^{(app)}$ should possess a jump
$J_{\lambda_{1}}^{(app)}$ across the circle $D_{\lambda_{1}}^{\epsilon}$. Besides, we obtain(see Appendix A)
\begin{align}\label{36}
J_{\lambda_{0}}^{(app)}-I=\left(\begin{array}{cc}
    O(t^{-\frac{1}{2}}) &  O(t^{-\frac{1}{2}-\mbox{Im}(\vartheta(\lambda_{0}))})\\
    O(t^{-\frac{1}{2}+\mbox{Im}(\vartheta(\lambda_{0}))})  &  O(t^{-\frac{1}{2}})\\
\end{array}\right),\notag\\
J_{\lambda_{1}}^{(app)}-I=\left(\begin{array}{cc}
    O(t^{-\frac{1}{2}}) &  O(t^{-\frac{1}{2}+\mbox{Im}(\vartheta(\lambda_{1}))})\\
    O(t^{-\frac{1}{2}-\mbox{Im}(\vartheta(\lambda_{1}))})  &  O(t^{-\frac{1}{2}})\\
\end{array}\right).
\end{align}
So the following RH problem  is given for matrix $M^{(app)(x, t, \lambda)}$
\begin{align}\label{37}
\left\{
\begin{array}{lr}
M^{(app)}(x, t, \lambda)\ \mbox{is analytic in} \ \mathbb{C }\setminus \Sigma^{(app)},\\
M_{+}^{(app)}(x, t, \lambda)=M_{-}^{(app)}(x, t, \lambda)J^{(app)}(x, t, \lambda), \qquad \lambda\in \Sigma^{(app)},\\
M^{(app)}(x, t, \lambda)\rightarrow I,\qquad \lambda\rightarrow \infty,\\
  \end{array}
\right.
\end{align}
of which the jump matrix $J^{(app)}(x, t, \lambda)$ is
\begin{align}\label{38}
J^{(app)}(x, t, \lambda)=\left\{
\begin{array}{lr}
J_{i}^{(app)}=J_{i}^{(2)}\ (i=1,2,3,4), \mbox{inside} \quad D_{\lambda_{0}}^{\epsilon},\\
J_{i}^{(app)}=J_{i}^{(2)}\ (i=5,6,7,8), \mbox{inside} \quad D_{\lambda_{1}}^{\epsilon},\\
J_{\lambda_{0}}^{(app)}\ \mbox{on} \quad D_{\lambda_{0}}^{\epsilon},\\
J_{\lambda_{1}}^{(app)}\  \mbox{on} \quad D_{\lambda_{1}}^{\epsilon},
  \end{array}
\right.
\end{align}
where
\begin{align}\label{39}
J_{1}^{(2)}=J_{5}^{(2)}=
\left(\begin{array}{cc}
    1 &  -\frac{r_{2}(\lambda)}{1-r_{1}(\lambda)r_{2}(\lambda)}e^{-2ift}\delta_{+}^{2}\\
    0  &  1\\
\end{array}\right),J_{2}^{(2)}=J_{6}^{(2)}=\left(\begin{array}{cc}
     1 &  0\\
     -r_{1}(\lambda)e^{2ift}\delta_{-}^{-2}  &  1\\
\end{array}\right), \notag\\
J_{3}^{(2)}=J_{7}^{(2)}=
\left(\begin{array}{cc}
    1 &  r_{2}(\lambda)e^{-2ift}\delta_{+}^{2}\\
    0  &  1\\
\end{array}\right),J_{4}^{(2)}=J_{8}^{(2)}=\left(\begin{array}{cc}
     1 &  0\\
     \frac{r_{1}(\lambda)}{1-r_{1}(\lambda)r_{2}(\lambda)}e^{2ift}\delta_{-}^{-2}  &  1\\
\end{array}\right).
\end{align}
\centerline{\begin{tikzpicture}[scale=2.0]
\draw[->][thick](0,1)--(0.5,0.5);
\draw[-][thick](0.5,0.5)--(1,0);
\draw[->][thick](0,-1)--(0.5,-0.5);
\draw[-][thick](0.5,-0.5)--(1,0);
\draw[->][thick](2,-1)--(1.5,-0.5);
\draw[-][thick](1.5,-0.5)--(1,0);
\draw[->][thick](2,1)--(1.5,0.5);
\draw[-][thick](1.5,0.5)--(1,0);
\draw[-][thick](0,1)--(-0.5,0.5);
\draw[->][thick](-1,0)--(-0.5,0.5);
\draw[-][thick](0,-1)--(-0.5,-0.5);
\draw[->][thick](-1,0)--(-0.5,-0.5);
\draw[->][thick](-1,0)--(-1.5,0.5);
\draw[-][thick](-1.5,0.5)--(-2,1);
\draw[->][thick](-1,0)--(-1.5,-0.5);
\draw[-][thick](-1.5,-0.5)--(-2,-1);
\draw[fill] (1,0)node[below]{$\lambda_{1}$} circle [radius=0.035];
\draw[fill] (-1,0)node[below]{$\lambda_{0}$} circle [radius=0.035];
\draw[-][thick](-0.5,0) arc(0:360:0.5);
\draw[-][thick](-0.5,0) arc(0:30:0.5);
\draw[-][thick](-0.5,0) arc(0:150:0.5);
\draw[-][thick](-0.5,0) arc(0:210:0.5);
\draw[-][thick](-0.5,0) arc(0:330:0.5);
\draw[-][thick](1.5,0) arc(0:360:0.5);
\draw[-][thick](1.5,0) arc(0:30:0.5);
\draw[-][thick](1.5,0) arc(0:150:0.5);
\draw[-][thick](1.5,0) arc(0:210:0.5);
\draw[-][thick](1.5,0) arc(0:330:0.5);
\draw[fill] (-0.8,0.25) node{$\textcolor[rgb]{1.00,0.00,0.00}{J_{1}^{(2)}}$};
\draw[fill] (-1.2,0.25) node{$\textcolor[rgb]{1.00,0.00,0.00}{J_{2}^{(2)}}$};
\draw[fill] (-1.25,-0.2) node{$\textcolor[rgb]{1.00,0.00,0.00}{J_{3}^{(2)}}$};
\draw[fill] (-0.7,-0.2) node{$\textcolor[rgb]{1.00,0.00,0.00}{J_{4}^{(2)}}$};
\draw[fill] (1.2,0.25) node{$\textcolor[rgb]{1.00,0.00,0.00}{J_{6}^{(2)}}$};
\draw[fill] (0.8,0.25) node{$\textcolor[rgb]{1.00,0.00,0.00}{J_{5}^{(2)}}$};
\draw[fill] (0.75,-0.2) node{$\textcolor[rgb]{1.00,0.00,0.00}{J_{8}^{(2)}}$};
\draw[fill] (1.3,-0.2) node{$\textcolor[rgb]{1.00,0.00,0.00}{J_{7}^{(2)}}$};
\end{tikzpicture}}
\centerline{\noindent { \textbf{Figure 4.} (Color online) The jump contour $\Sigma^{(2)}$.}}

For large $\lambda$, we define the factorization
\begin{align}\label{40}
M^{(2)}=M^{(err)}M^{(app)},
\end{align}
where the error term contains higher-order contribution from the contour $\Sigma^{(2)}$. Then  matrix $M^{(err)}(x, t, \lambda)$ meets the following RH problem:
\begin{align}\label{41}
\left\{
\begin{array}{lr}
M^{(err)}(x, t, \lambda)\ \mbox{is analytic in} \ \mathbb{C }\setminus \Sigma^{(err)},\\
M_{+}^{(err)}(x, t, \lambda)=M_{-}^{(err)}(x, t, \lambda)J^{(err)}(x, t, \lambda), \qquad \lambda\in \Sigma^{(err)},\\
M^{(err)}(x, t, \lambda)\rightarrow I,\qquad \lambda\rightarrow \infty,\\
  \end{array}
\right.
\end{align}
of which the jump matrix $J^{(err)}(x, t, \lambda)$ is(see Appendix B)
\begin{align}\label{42}
J^{(err)}(x, t, \lambda)=\left\{
\begin{array}{lr}
J_{i}^{(err)}=J_{i}^{(2)}=I+O(e^{-\tilde{C}t}),\ (i=1,2,\cdots, 8), \mbox{outside} \ D_{\lambda_{0}}^{\epsilon}\cup D_{\lambda_{1}}^{\epsilon},\\
J_{\lambda_{0}}^{(err)}=(J_{\lambda_{0}}^{(app)})^{-1}\ \mbox{on} \quad D_{\lambda_{0}}^{\epsilon},\\
J_{\lambda_{1}}^{(err)}=(J_{\lambda_{1}}^{(app)})^{-1}\  \mbox{on} \quad D_{\lambda_{1}}^{\epsilon}.
  \end{array}
\right.
\end{align}
Let's expand the matrices $M^{(2)}, M_{\lambda_{0}}^{(app)}, M_{\lambda_{1}}^{(app)}, M^{(err)}$ at infinity into the Laurent series
\begin{align}\label{43}
&M^{(2)}=I+\frac{M_{1}^{(2)}}{\lambda}+\frac{M_{2}^{(2)}}{\lambda^{2}}+\cdots, \qquad \lambda\rightarrow\infty,\notag\\
&M_{\lambda_{0}}^{(app)}=I+\frac{(M_{\lambda_{0}}^{(app)})_{1}}{\lambda}+\frac{(M_{\lambda_{0}}^{(app)})_{2}}{\lambda^{2}}+\cdots, \qquad \lambda\rightarrow\infty,\notag\\
&M_{\lambda_{1}}^{(app)}=I+\frac{(M_{\lambda_{1}}^{(app)})_{1}}{\lambda}+\frac{(M_{\lambda_{1}}^{(app)})_{2}}{\lambda^{2}}+\cdots, \qquad \lambda\rightarrow\infty,\notag\\
&M^{(err)}=I+\frac{M_{1}^{(err)}}{\lambda}+\frac{M_{2}^{(err)}}{\lambda^{2}}+\cdots, \qquad \lambda\rightarrow\infty.
\end{align}
According to the factorization \eqref{40}, comparing the coefficients of $\frac{1}{\lambda}$, we find
\begin{align}\label{44}
M_{1}^{(2)}=(M_{\lambda_{0}}^{(app)})_{1}+(M_{\lambda_{1}}^{(app)})_{1}+M_{1}^{(err)}.
\end{align}
Thus the solution of the nonlocal Hirota equation \eqref{1} is
\begin{align}\label{45}
q(x,t)=2i\left[M_{1}^{(2)}\right]_{12}=2i\left[(M_{\lambda_{0}}^{(app)})_{1}\right]_{12}
+2i\left[(M_{\lambda_{1}}^{(app)})_{1}\right]_{12}+2i\left[M_{1}^{(err)}\right]_{12}.
\end{align}
Similar to literature \cite{Wangdengshan-na24,Peng16-19},
the absolution of matrix $M_{1}^{(err)}(x, t, \lambda)$ in Eq.\eqref{41} satisfies(see Appendix C)
\begin{align}\label{46}
|M_{1}^{(err)}(x, t, \lambda)|=
O(t^{-\frac{1}{2}-\mbox{max}\{|\mbox{Im}\vartheta(\lambda_{0})|,|\mbox{Im}\vartheta(\lambda_{1})|\}}).
\end{align}

\subsection{Reduction to a model Riemann-Hilbert problem}
In this subsection, we will define a scaling transformation to separate the time $t$ from the
jump matrix, given by
\begin{align}\label{47}
\tilde{\lambda}=T_{0}(\lambda)=\sqrt{-8t(\alpha+6\beta \lambda_{0})}(\lambda-\lambda_{0}), \quad \lambda\in \Sigma_{\lambda_{0}}^{(app)},\notag\\
\tilde{\lambda}=T_{1}(\lambda)=\sqrt{8t(\alpha+6\beta \lambda_{1})}(\lambda-\lambda_{1})\quad \lambda\in \Sigma_{\lambda_{1}}^{(app)}.
\end{align}
For a given  function $\varphi(\zeta)$, one has
$T_{0}(\varphi(\lambda))=\varphi(\frac{\tilde{\lambda}}{\sqrt{-8t(\alpha+6\beta \lambda_{0})}}+\lambda_{0})$ and $T_{1}(\varphi(\lambda))=\varphi(\frac{\tilde{\lambda}}{\sqrt{8t(\alpha+6\beta \lambda_{1})}}+\lambda_{1})$.
Hence, we have
\begin{align}
T_{0}(e^{-itf}\delta(\lambda))
=\delta_{\lambda_{0}}^{0}\delta_{\lambda_{0}}^{1},\quad
T_{1}(e^{-itf}\delta(\lambda))
=\delta_{\lambda_{1}}^{0}\delta_{\lambda_{1}}^{1}.
\end{align}
where
\begin{align}\label{48}
\delta_{\lambda_{0}}^{0}&=\left[-32\lambda_{0}^{2}t(\alpha+6\beta\lambda_{0})\right]^{\frac{i\vartheta(\lambda_{0})}{2}}
e^{2i\lambda_{0}^{2}t(4\beta\lambda_{0}+\alpha)+\chi_{0}(\lambda_{0})},\notag\\
\delta_{\lambda_{0}}^{1}&=\tilde{\lambda}^{-i\vartheta(\lambda_{0})}\left(\frac{2\lambda_{0}}{\tilde{\lambda}/\sqrt{-8t(\alpha+6\beta\lambda_{0})}+\lambda_{0}-\lambda_{1}}\right)
^{-i\vartheta(\lambda_{0})}\notag\\
&e^{\frac{i}{4}\tilde{\lambda}^{2}\left(1-\frac{i\beta \tilde{\lambda}}{\sqrt{2t}(6\beta\lambda_{0}+\alpha)^{\frac{3}{2}}}\right)}
e^{\chi_{0}(\frac{\tilde{\lambda}}{\sqrt{-8t(\alpha+6\beta \lambda_{0})}}+\lambda_{0})-\chi_{0}(\lambda_{0})},\notag\\
\delta_{\lambda_{1}}^{0}&=\left[32\lambda_{1}^{2}t(\alpha+6\beta\lambda_{1})\right]^{-\frac{i\vartheta(\lambda_{1})}{2}}
e^{2i\lambda_{1}^{2}t(4\beta\lambda_{1}+\alpha)+\chi_{1}(\lambda_{1})},\notag\\
\delta_{\lambda_{1}}^{1}&=\tilde{\lambda}^{i\vartheta(\lambda_{1})}\left(\frac{2\lambda_{1}}{\tilde{\lambda}/\sqrt{8t(\alpha+6\beta\lambda_{1})}+\lambda_{1}-\lambda_{0}}\right)
^{i\vartheta(\lambda_{1})}\notag\\
&e^{-\frac{i}{4}\tilde{\lambda}^{2}\left(1+\frac{\beta \tilde{\lambda}}{\sqrt{2t}(6\beta\lambda_{1}+\alpha)^{\frac{3}{2}}}\right)}
e^{\chi_{1}(\frac{\tilde{\lambda}}{\sqrt{8t(\alpha+6\beta \lambda_{1})}}+\lambda_{1})-\chi_{1}(\lambda_{1})}.
\end{align}

\subsubsection{The $T_{0}$ scaling transformation}
Here, we first consider the  scaling transformation of $T_{0}$, it is easy to obtain the following RH problem
\begin{align}\label{49}
M^{(3)}_{+}(x, t; \tilde{\lambda})=M^{(3)}_{-}(x, t; \tilde{\lambda})J^{(3)}, \quad \tilde{\lambda}\in\Sigma^{(3)},
\end{align}
where we have defined
\begin{align}
M^{(3)}(x, t; \tilde{\lambda})=T_{0}(M_{\lambda_{0}}^{(app)}(x, t; \lambda)),\qquad J^{(3)}(x, t; \tilde{\lambda})=T_{0}(J_{i}^{(2)}(x, t; \lambda)),\quad i=1,2,3,4.\notag
\end{align}
In terms of the above analysis, one gets the new jump matrix $J^{(3)}$, given by(see Figure 5)
\begin{align}\label{50}
&J_{1}^{(3)}=
\left(\begin{array}{cc}
    1  &  -(\delta_{\lambda_{0}}^{0}\delta_{\lambda_{0}}^{1})^{2}\frac{r_{2}}{1-r_{1}r_{2}}(\frac{\tilde{\lambda}}{\sqrt{-8t(\alpha+6\beta \lambda_{0})}}+\lambda_{0})\\
    0  &  1\\
\end{array}\right),\notag\\
&J_{2}^{(3)}=
\left(\begin{array}{cc}
    1  &  0\\
    -(\delta_{\lambda_{0}}^{0}\delta_{\lambda_{0}}^{1})^{-2}r_{1}(\frac{\tilde{\lambda}}{\sqrt{-8t(\alpha+6\beta\lambda_{0})}}+\lambda_{0})  &  1\\
\end{array}\right), \notag\\
&J_{3}^{(3)}=\left(\begin{array}{cc}
    1  &  (\delta_{\lambda_{0}}^{0}\delta_{\lambda_{0}}^{1})^{2}r_{2}(\frac{\tilde{\lambda}}{\sqrt{-8t(\alpha+6\beta\lambda_{0})}}+\lambda_{0})\\
    0  &  1\\
\end{array}\right),\notag\\
&J_{4}^{(3)}=\left(\begin{array}{cc}
    1  &  0\\
    (\delta_{\lambda_{0}}^{0}\delta_{\lambda_{0}}^{1})^{-2}\frac{r_{1}}{1-r_{1}r_{2}}(\frac{\tilde{\lambda}}{\sqrt{-8t(\alpha+6\beta \lambda_{0})}}+\lambda_{0})  &  1\\
\end{array}\right).
\end{align}

\centerline{\begin{tikzpicture}
\draw[fill] (0,0) circle [radius=0.035];
\draw[-][thick](-2.12,-2.12)node[left]{$J_{3}^{(3)}$}--(-1.06,-1.06);
\draw[<-][thick](-1.06,-1.06)--(0,0)node[below]{$O$};
\draw[->][thick](0,0)--(1.06,1.06);
\draw[-][thick](1.06,1.06)--(2.12,2.12)node[right]{$J_{1}^{(3)}$};
\draw[-][thick](-2.12,2.12)node[left]{$J_{2}^{(3)}$}--(-1.06,1.06);
\draw[<-][thick](-1.06,1.06)--(0,0);
\draw[->][thick](0,0)--(1.06,-1.06);
\draw[-][thick](1.06,-1.06)--(2.12,-2.12)node[right]{$J_{4}^{(3)}$};
\end{tikzpicture}}

\centerline{\noindent {\small \textbf{Figure 5.} (Color online) The jump contour $\Sigma^{(3)}$.}}

Since
\begin{align}
M^{(3)}&=T_{0}(M_{\lambda_{0}}^{(app)}(\lambda))=M_{\lambda_{0}}^{(app)}(\frac{\tilde{\lambda}}{\sqrt{-8t(\alpha+6\beta \lambda_{0})}}+\lambda_{0})\notag\\
&=I+\frac{(M_{\lambda_{0}}^{(app)})_{1}}{\frac{\tilde{\lambda}}{\sqrt{-8t(\alpha+6\beta \lambda_{0})}}+\lambda_{0}}+\cdots=I+\frac{M_{1}^{(3)}}{\tilde{\lambda}}+\cdots,
\end{align}
Comparing the coefficient of $\tilde{\lambda}$ in above formulas, we have
\begin{align}\label{50.1}
M_{1}^{(3)}=\sqrt{-8t(\alpha+6\beta\lambda_{0})}(M_{\lambda_{0}}^{(app)})_{1}.
\end{align}

Moreover, as $t\rightarrow\infty$, one obtains
\begin{align}
\lim_{t\rightarrow \infty}(\frac{\tilde{\lambda}}{\sqrt{-8t(\alpha+6\beta\lambda_{0})}}+\lambda_{0})=\lambda_{0},\notag\\
\lim_{t\rightarrow \infty}\delta_{\lambda_{0}}^{1}=\tilde{\lambda}^{-i\vartheta(\lambda_{0})}e^{\frac{1}{4}i\tilde{\lambda}^{2}}.
\end{align}
To separate the time $t$ completely, we perform the following limiting operation
\begin{align}\label{50.2}
M^{(\infty)}=\lim_{t\rightarrow\infty}(\delta_{\lambda_{0}}^{0})^{-\hat{\sigma}_{3}}M^{(3)},
\end{align}
which changes the jumping curve $\Sigma^{(3)}$ into $\Sigma^{(\infty)}$, and leads to the following RH problem:\\
\begin{align}\label{51}
\left\{
\begin{array}{lr}
M^{(\infty)}(x, t; \tilde{\lambda})\ \mbox{is analytic in} \ \mathbb{C}\backslash\Sigma^{(\infty)},\\
M^{(\infty)}_{+}(x, t; \tilde{\lambda})=M^{(\infty)}_{-}(x, t; \tilde{\lambda})J^{(\infty)}(x, t; \tilde{\lambda}), \quad \tilde{\lambda}\in\Sigma^{(\infty)},\\
M^{(\infty)}(x, t; \tilde{\lambda})\rightarrow I,\qquad \tilde{\lambda}\rightarrow \infty,
  \end{array}
\right.
\end{align}
where
\begin{align}
&J_{1}^{(\infty)}=
\left(\begin{array}{cc}
    1  &  -\tilde{\lambda}^{-2i\vartheta(\lambda_{0})}e^{\frac{1}{2}i\tilde{\lambda}^{2}}\frac{r_{2}}{1-r_{1}r_{2}}(\lambda_{0})\\
    0  &  1\\
\end{array}\right),\notag\\
&J_{2}^{(\infty)}=
\left(\begin{array}{cc}
    1  &  0\\
    -\tilde{\lambda}^{2i\vartheta(\lambda_{0})}e^{-\frac{1}{2}i\tilde{\lambda}^{2}}r_{1}(\lambda_{0})  &  1\\
\end{array}\right), \notag\\
&J_{3}^{(\infty)}=\left(\begin{array}{cc}
    1  &  \tilde{\lambda}^{-2i\vartheta(\lambda_{0})}e^{\frac{1}{2}i\tilde{\lambda}^{2}}r_{2}(\lambda_{0})\\
    0  &  1\\
\end{array}\right),\notag\\
&J_{4}^{(\infty)}=\left(\begin{array}{cc}
    1  &  0\\
    \tilde{\lambda}^{2i\vartheta(\lambda_{0})}e^{-\frac{1}{2}i\tilde{\lambda}^{2}}\frac{r_{1}}{1-r_{1}r_{2}}(\lambda_{0})  &  1\\
\end{array}\right).
\end{align}

To obtain the model RH problem, we define the following transformation
\begin{align}\label{52}
M^{(mod)}=M^{(\infty)}G_{j},\qquad  \tilde{\lambda} \in \Omega_{j},\qquad j=0, \cdots, 4,
\end{align}
where
\begin{gather}\label{53}
G_{0}=e^{\frac{1}{4}i\tilde{\lambda}^{2}\sigma_{3}}\tilde{\lambda}^{-i\vartheta(\lambda_{0})\sigma_{3}},\notag\\
G_{1}=G_{0}\left(\begin{array}{cc}
    1  &  -\frac{r_{2}}{1-r_{1}r_{2}}(\lambda_{0})\\
    0  &  1\\
\end{array}\right),\quad G_{2}=G_{0}\left(\begin{array}{cc}
    1  &  0\\
    r_{1}(\lambda_{0})  &  1\\
\end{array}\right),\notag\\
G_{3}=G_{0}\left(\begin{array}{cc}
    1  &  r_{2}(\lambda_{0})\\
    0  &  1\\
\end{array}\right),\quad G_{4}=G_{0}\left(\begin{array}{cc}
    1  &  0\\
    -\frac{r_{1}}{1-r_{1}r_{2}}(\lambda_{0})  &  1\\
\end{array}\right).
\end{gather}
Through this transformation, we obtain a model RH problem for $M^{(mod)}$ with a constant jump matrix, given by
\begin{align}\label{54}
\left\{
\begin{array}{lr}
M^{(mod)}(x, t; \tilde{\lambda})\ \mbox{is analytic in} \ \mathbb{C}\backslash \mathbb{R},\\
M^{(mod)}_{+}(x, t; \tilde{\lambda})=M^{(mod)}_{-}(x, t; \tilde{\lambda})J^{(mod)}(x, t; \tilde{\lambda}), \quad \tilde{\lambda}\in \mathbb{R},\\
M^{(mod)}(x, t; \tilde{\lambda})\rightarrow e^{\frac{1}{4}i\tilde{\lambda}^{2}\sigma_{3}}\tilde{\lambda}^{-i\vartheta(\lambda_{0})\sigma_{3}},\qquad \tilde{\lambda}\rightarrow \infty,
  \end{array}
\right.
\end{align}
where
\begin{align}
J^{(mod)}=\left(\begin{array}{cc}
    1-r_{1}(\lambda_{0})r_{2}(\lambda_{0})  &  -r_{2}(\lambda_{0})\\
  r_{1}(\lambda_{0}) &  1\\
\end{array}\right).
\end{align}
The solution $M^{(mod)}(\tilde{\lambda})$ of this RH problem can be given explicitly via using the parabolic cylinder functions.

In order to  derive the asymptotic formulas in Theorem 1, we give the large-$\tilde{\lambda}$ behavior
of $M^{(\infty)}(\tilde{\lambda})$(see  Appendix D)
\begin{align}\label{55}
M^{(\infty)}(\tilde{\lambda})=I+\frac{M_{1}^{(\infty)}}{\tilde{\lambda}}+O(\frac{1}{\tilde{\lambda}}),\qquad \tilde{\lambda}\rightarrow\infty,
\end{align}
where
\begin{align}\label{56}
&\left[M_{1}^{(\infty)}\right]_{12}=\frac{\sqrt{2\pi}i(-1)^{i\vartheta(\lambda_{0})-\frac{1}{2}}e^{-\frac{3\pi i}{4}+\frac{\pi\vartheta(\lambda_{0}) }{2}}}{r_{1}(\lambda_{0})\Gamma(i\vartheta(\lambda_{0}))},\notag\\ &\left[M_{1}^{(\infty)}\right]_{21}=\frac{\sqrt{2\pi}i(-1)^{-i\vartheta(\lambda_{0})-\frac{1}{2}}e^{-\frac{\pi i}{4}+\frac{\pi\vartheta(\lambda_{0}) }{2}}}{r_{2}(\lambda_{0})\Gamma(-i\vartheta(\lambda_{0}))}.
\end{align}
Therefore, combining \eqref{50.1}, \eqref{50.2}, \eqref{55} and \eqref{56}, we have
\begin{align}\label{57}
&\left[(M_{\lambda_{0}}^{(app)})_{1}\right]_{12}=\frac{1}{\sqrt{-8t(\alpha+6\beta\lambda_{0})}}\left[M_{1}^{(3)}\right]_{12}
=\frac{(\delta_{\lambda_{0}}^{0})^{2}}{\sqrt{-8t(\alpha+6\beta\lambda_{0})}}\left[M_{1}^{(\infty)}\right]_{12}\notag\\
&=\frac{\sqrt{2\pi}\left[32\lambda_{0}^{2}t(\alpha+6\beta\lambda_{0})\right]^{i\vartheta(\lambda_{0})}
e^{4i\lambda_{0}^{2}t(4\beta\lambda_{0}+\alpha)+2\chi_{0}(\lambda_{0})-\frac{3\pi i}{4}+\frac{\pi\vartheta(\lambda_{0}) }{2}}}{\sqrt{-8t(\alpha+6\beta\lambda_{0})}r_{1}(\lambda_{0})\Gamma(i\vartheta(\lambda_{0}))}.
\end{align}

\centerline{\begin{tikzpicture}
\draw[fill] (0,0) circle [radius=0.035];
\draw[-][thick](-2.12,-2.12)--(-1.06,-1.06);
\draw[<-][thick](-1.06,-1.06)--(0,0)node[below]{$O$};
\draw[->][thick](0,0)--(1.06,1.06);
\draw[-][thick](1.06,1.06)--(2.12,2.12);
\draw[-][thick](-2.12,2.12)--(-1.06,1.06);
\draw[<-][thick](-1.06,1.06)--(0,0);
\draw[->][thick](0,0)--(1.06,-1.06);
\draw[->][thick](-3,0)--(-1.5,0);
\draw[->][thick](-1.5,0)--(1.5,0);
\draw[-][thick](1.5,0)--(3.0,0)node[right]{$\mathbb{R}$};
\draw[-][thick](1.06,-1.06)--(2.12,-2.12);
\draw[fill] (0,1) node{$\Omega_{0}$};
\draw[fill] (0,-1) node{$\Omega_{0}$};
\draw[fill] (-1,0.5) node{$\Omega_{2}$};
\draw[fill] (1,0.5) node{$\Omega_{1}$};
\draw[fill] (-1,-0.5) node{$\Omega_{3}$};
\draw[fill] (1,-0.5) node{$\Omega_{4}$};
\end{tikzpicture}}
\centerline{\noindent {\small \textbf{Figure 6.} (Color online) The jump contour $\Sigma^{(3)}\cup\Sigma^{(mod)}$.}}

\subsubsection{The $T_{1}$ scaling transformation}
Here, we perform the $T_{1}$ scaling transformation, and obtain the following RH problem
\begin{align}\label{58}
M^{(3)}_{+}(x, t; \tilde{\lambda})=M^{(3)}_{-}(x, t; \tilde{\lambda})J^{(3)}, \quad \tilde{\lambda}\in\Sigma^{(3)},
\end{align}
where
\begin{align}
M^{(3)}(x, t; \tilde{\lambda})=T_{1}(M_{\lambda_{1}}^{(app)}(x, t; \lambda)),\qquad J^{(3)}(x, t; \tilde{\lambda})=T_{1}(J_{i}^{(2)}(x, t; \lambda)),\quad i=5,6,7,8.\notag
\end{align}
and the new jump matrix $J^{(3)}$ is (see Figure 7)
\begin{align}\label{59}
&J_{5}^{(3)}=
\left(\begin{array}{cc}
    1  &  -(\delta_{\lambda_{1}}^{0}\delta_{\lambda_{1}}^{1})^{2}\frac{r_{2}}{1-r_{1}r_{2}}(\frac{\tilde{\lambda}}{\sqrt{8t(\alpha+6\beta \lambda_{1})}}+\lambda_{1})\\
    0  &  1\\
\end{array}\right),\notag\\
&J_{6}^{(3)}=
\left(\begin{array}{cc}
    1  &  0\\
    -(\delta_{\lambda_{1}}^{0}\delta_{\lambda_{1}}^{1})^{-2}r_{1}(\frac{\tilde{\lambda}}{\sqrt{8t(\alpha+6\beta \lambda_{1})}}+\lambda_{1})  &  1\\
\end{array}\right), \notag\\
&J_{7}^{(3)}=\left(\begin{array}{cc}
    1  &  (\delta_{\lambda_{1}}^{0}\delta_{\lambda_{1}}^{1})^{2}r_{2}(\frac{\tilde{\lambda}}{\sqrt{8t(\alpha+6\beta \lambda_{1})}}+\lambda_{1})\\
    0  &  1\\
\end{array}\right),\notag\\
&J_{8}^{(3)}=\left(\begin{array}{cc}
    1  &  0\\
    (\delta_{\lambda_{1}}^{0}\delta_{\lambda_{1}}^{1})^{-2}\frac{r_{1}}{1-r_{1}r_{2}}(\frac{\tilde{\lambda}}{\sqrt{8t(\alpha+6\beta \lambda_{1})}}+\lambda_{1})  &  1\\
\end{array}\right).
\end{align}

\centerline{\begin{tikzpicture}
\draw[fill] (0,0) circle [radius=0.035];
\draw[->][thick](-2.12,-2.12)node[left]{$J_{8}^{(3)}$}--(-1.06,-1.06);
\draw[-][thick](-1.06,-1.06)--(0,0)node[below]{$O$};
\draw[-][thick](0,0)--(1.06,1.06);
\draw[<-][thick](1.06,1.06)--(2.12,2.12)node[right]{$J_{6}^{(3)}$};
\draw[->][thick](-2.12,2.12)node[left]{$J_{5}^{(3)}$}--(-1.06,1.06);
\draw[-][thick](-1.06,1.06)--(0,0);
\draw[-][thick](0,0)--(1.06,-1.06);
\draw[<-][thick](1.06,-1.06)--(2.12,-2.12)node[right]{$J_{7}^{(3)}$};
\end{tikzpicture}}

\centerline{\noindent {\small \textbf{Figure 7.} (Color online) The jump contour $\Sigma^{(3)}$.}}

Due to
\begin{align}
M^{(3)}&=T_{1}(M_{\lambda_{1}}^{(app)}(\lambda))=M_{\lambda_{1}}^{(app)}(\frac{\tilde{\lambda}}{\sqrt{8t(\alpha+6\beta \lambda_{1})}}+\lambda_{1})\notag\\
&=I+\frac{(M_{\lambda_{1}}^{(app)})_{1}}{\frac{\tilde{\lambda}}{\sqrt{8t(\alpha+6\beta \lambda_{1})}}+\lambda_{1}}+\cdots=I+\frac{M_{1}^{(3)}}{\tilde{\lambda}}+\cdots,
\end{align}
we have
\begin{align}\label{60}
M_{1}^{(3)}=\sqrt{8t(\alpha+6\beta \lambda_{1})}(M_{\lambda_{1}}^{(app)})_{1}.
\end{align}
Besides, we arrive at
\begin{align}
&\lim_{t\rightarrow \infty}(\frac{\tilde{\lambda}}{\sqrt{8t(\alpha+6\beta \lambda_{1})}}+\lambda_{1})=\lambda_{1},\notag\\
&\lim_{t\rightarrow \infty}\delta_{\lambda_{1}}^{1}=\tilde{\lambda}^{i\vartheta(\lambda_{1})}e^{-\frac{1}{4}i\tilde{\lambda}^{2}}.
\end{align}
Similarly, we carry out the following limitation
\begin{align}\label{61}
M^{(\infty)}=\lim_{t\rightarrow\infty}(\delta_{\lambda_{1}}^{0})^{-\hat{\sigma}_{3}}M^{(3)},
\end{align}
and we obtain the following new RH problem:\\
\begin{align}\label{62}
\left\{
\begin{array}{lr}
M^{(\infty)}(x, t; \tilde{\lambda})\ \mbox{is analytic in} \ \mathbb{C}\backslash\Sigma^{(\infty)},\\
M^{(\infty)}_{+}(x, t; \tilde{\lambda})=M^{(\infty)}_{-}(x, t; \tilde{\lambda})J^{(\infty)}(x, t; \tilde{\lambda}), \quad \tilde{\lambda}\in\Sigma^{(\infty)},\\
M^{(\infty)}(x, t; \tilde{\lambda})\rightarrow I,\qquad \tilde{\lambda}\rightarrow \infty,
  \end{array}
\right.
\end{align}
where
\begin{align}
&J_{5}^{(\infty)}=
\left(\begin{array}{cc}
    1  &  -\tilde{\lambda}^{2i\vartheta(\lambda_{1})}e^{-\frac{1}{2}i\tilde{\lambda}^{2}}\frac{r_{2}}{1-r_{1}r_{2}}(\lambda_{1})\\
    0  &  1\\
\end{array}\right),
J_{6}^{(\infty)}=
\left(\begin{array}{cc}
    1  &  0\\
    -\tilde{\lambda}^{-2i\vartheta(\lambda_{1})}e^{\frac{1}{2}i\tilde{\lambda}^{2}}r_{1}(\lambda_{1})  &  1\\
\end{array}\right), \notag\\
&J_{7}^{(\infty)}=\left(\begin{array}{cc}
    1  &  \tilde{\lambda}^{2i\vartheta(\lambda_{1})}e^{-\frac{1}{2}i\tilde{\lambda}^{2}}r_{2}(\lambda_{1})\\
    0  &  1\\
\end{array}\right),
J_{8}^{(\infty)}=\left(\begin{array}{cc}
    1  &  0\\
    \tilde{\lambda}^{-2i\vartheta(\lambda_{1})}e^{\frac{1}{2}i\tilde{\lambda}^{2}}\frac{r_{1}}{1-r_{1}r_{2}}(\lambda_{1})  &  1\\
\end{array}\right).
\end{align}
Next, we will aim to obtain a model RH problem by defining
\begin{align}\label{63}
M^{(mod)}=M^{(\infty)}G_{j},\qquad  \tilde{\lambda} \in \Omega_{j},\qquad j=0, \cdots, 4,
\end{align}
where
\begin{gather}\label{64}
G_{0}=e^{-\frac{1}{4}i\tilde{\lambda}^{2}\sigma_{3}}\tilde{\lambda}^{i\vartheta(\lambda_{1})\sigma_{3}},\notag\\
G_{1}=G_{0}\left(\begin{array}{cc}
    1  &  -\frac{r_{2}}{1-r_{1}r_{2}}(\lambda_{1})\\
    0  &  1\\
\end{array}\right),\quad G_{2}=G_{0}\left(\begin{array}{cc}
    1  &  0\\
    r_{1}(\lambda_{1})  &  1\\
\end{array}\right),\notag\\
G_{3}=G_{0}\left(\begin{array}{cc}
    1  &  r_{2}(\lambda_{1})\\
    0  &  1\\
\end{array}\right),\quad G_{4}=G_{0}\left(\begin{array}{cc}
    1  &  0\\
    -\frac{r_{1}}{1-r_{1}r_{2}}(\lambda_{1})  &  1\\
\end{array}\right),
\end{gather}
then we obtain the following model RH problem:
\begin{align}\label{65}
\left\{
\begin{array}{lr}
M^{(mod)}(x, t; \tilde{\lambda})\ \mbox{is analytic in} \ \mathbb{C}\backslash \mathbb{R},\\
M^{(mod)}_{+}(x, t; \tilde{\lambda})=M^{(mod)}_{-}(x, t; \tilde{\lambda})J^{(mod)}(x, t; \tilde{\lambda}), \quad \tilde{\lambda}\in \mathbb{R},\\
M^{(mod)}(x, t; \tilde{\lambda})\rightarrow e^{-\frac{1}{4}i\tilde{\lambda}^{2}\sigma_{3}}\tilde{\lambda}^{i\vartheta(\lambda_{1})\sigma_{3}},\qquad \tilde{\lambda}\rightarrow \infty,
  \end{array}
\right.
\end{align}
where
\begin{align}
J^{(mod)}=\left(\begin{array}{cc}
    1-r_{1}(\lambda_{1})r_{2}(\lambda_{1})  &  -r_{2}(\lambda_{1})\\
  r_{1}(\lambda_{1}) &  1\\
\end{array}\right).
\end{align}
Performing the same procedure in Appendix D, we get the large-$\tilde{\lambda}$ behavior
of $M^{(\infty)}(\tilde{\lambda})$:
\begin{align}\label{66}
M^{(\infty)}(\tilde{\lambda})=I+\frac{M_{1}^{(\infty)}}{\tilde{\lambda}}+O(\frac{1}{\tilde{\lambda}}),\qquad \tilde{\lambda}\rightarrow\infty,
\end{align}
where
\begin{align}\label{67}
\left[M_{1}^{(\infty)}\right]_{12}=\frac{\sqrt{2\pi}ie^{-\frac{\pi}{2}\vartheta(\lambda_{1})}e^{-\frac{3\pi i}{4}}}{r_{1}(\lambda_{1})\Gamma(-i\vartheta(\lambda_{1}))},\quad \left[M_{1}^{(\infty)}\right]_{21}=\frac{\sqrt{2\pi}ie^{-\frac{\pi}{2}\vartheta(\lambda_{1})}e^{-\frac{\pi i}{4}}}{r_{2}(\lambda_{1})\Gamma(i\vartheta(\lambda_{1}))}.
\end{align}
Therefore, combining \eqref{60}, \eqref{61}, \eqref{66} and \eqref{67}, we have
\begin{align}\label{68}
&\left[(M_{\lambda_{1}}^{(app)})_{1}\right]_{12}=\frac{1}{\sqrt{8t(\alpha+6\beta\lambda_{1})}}\left[M_{1}^{(3)}\right]_{12}
=\frac{(\delta_{\lambda_{1}}^{0})^{2}}{\sqrt{8t(\alpha+6\beta\lambda_{1})}}\left[M_{1}^{(\infty)}\right]_{12}\notag\\
&=\frac{\sqrt{2\pi}i\left[32\lambda_{1}^{2}t(\alpha+6\beta\lambda_{1})\right]^{-i\vartheta(\lambda_{1})}
e^{4i\lambda_{1}^{2}t(4\beta\lambda_{1}+\alpha)+2\chi_{1}(\lambda_{1})-\frac{\pi}{2}\vartheta(\lambda_{1})-\frac{3\pi i}{4} }}{\sqrt{8t(\alpha+6\beta\lambda_{1})}r_{1}(\lambda_{1})\Gamma(-i\vartheta(\lambda_{1}))}.
\end{align}
Finally, combining  \eqref{45}, \eqref{46}, \eqref{57} and \eqref{68},  we can achieve the result of Theorem 1.1.

\centerline{\begin{tikzpicture}
\draw[fill] (0,0) circle [radius=0.035];
\draw[->][thick](-2.12,-2.12)--(-1.06,-1.06);
\draw[-][thick](-1.06,-1.06)--(0,0)node[below]{$O$};
\draw[-][thick](0,0)--(1.06,1.06);
\draw[<-][thick](1.06,1.06)--(2.12,2.12);
\draw[->][thick](-2.12,2.12)--(-1.06,1.06);
\draw[-][thick](-1.06,1.06)--(0,0);
\draw[-][thick](0,0)--(1.06,-1.06);
\draw[->][thick](-3,0)--(-1.5,0);
\draw[->][thick](-1.5,0)--(1.5,0);
\draw[-][thick](1.5,0)--(3.0,0)node[right]{$\mathbb{R}$};
\draw[<-][thick](1.06,-1.06)--(2.12,-2.12);
\draw[fill] (0,1) node{$\Omega_{0}$};
\draw[fill] (0,-1) node{$\Omega_{0}$};
\draw[fill] (-1,0.5) node{$\Omega_{1}$};
\draw[fill] (1,0.5) node{$\Omega_{2}$};
\draw[fill] (-1,-0.5) node{$\Omega_{4}$};
\draw[fill] (1,-0.5) node{$\Omega_{3}$};
\end{tikzpicture}}
\centerline{\noindent {\small \textbf{Figure 8.} (Color online) The jump contour $\Sigma^{(3)}\cup\Sigma^{(mod)}$.}}

\section*{Appendix A}
From Ref.\cite{DeiftZ}, we have
\begin{align}\label{75}
M^{(\infty)}=(\delta_{\lambda_{0}}^{0})^{-\hat{\sigma}_{3}}M^{(3)}+O(t^{-\frac{1}{2}}\ln t).
\end{align}
For $\tilde{\lambda}\rightarrow \infty$
\begin{align}\label{76}
M^{(3)}=I+M_{1}^{(3)}\tilde{\lambda}^{-1}+O(\tilde{\lambda}^{-2})
=I+(\delta_{\lambda_{0}}^{0})^{\hat{\sigma}_{3}}(M_{1}^{(\infty)}+O(t^{-\frac{1}{2}}\ln t))\tilde{\lambda}^{-1}+O(\tilde{\lambda}^{-2}).
\end{align}
On $D_{\lambda_{0}}^{\epsilon}$, $M_{\lambda_{0}-}^{(app)}=I$, one has
\begin{align}\label{77}
M_{\lambda_{0}+}^{(app)}&=T_{0}^{-1}(M^{(3)})
=I+(\delta_{\lambda_{0}}^{0})^{\hat{\sigma}_{3}}(T_{0}^{-1}(M_{1}^{(\infty)})+O(t^{-\frac{1}{2}}\ln t))(\sqrt{8t(\alpha+6\beta |\lambda_{0}|)}(\lambda-\lambda_{0}))^{-1}\notag\\
&+O((\sqrt{8t(\alpha+6\beta |\lambda_{0}|)}(\lambda-\lambda_{0}))^{-2}).
\end{align}
Therefore, we get
\begin{align}\label{78}
J_{\lambda_{0}}^{(app)}-I=(M_{\lambda_{0}-}^{(app)})^{-1}M_{\lambda_{0}+}^{(app)}-I=\left(\begin{array}{cc}
    O(t^{-\frac{1}{2}}) &  O(t^{-\frac{1}{2}-\mbox{Im}(\vartheta(\lambda_{0}))})\\
    O(t^{-\frac{1}{2}+\mbox{Im}(\vartheta(\lambda_{0}))})  &  O(t^{-\frac{1}{2}})\\
\end{array}\right).
\end{align}
In a similar way, we have
\begin{align}\label{79}
J_{\lambda_{1}}^{(app)}-I=\left(\begin{array}{cc}
    O(t^{-\frac{1}{2}}) &  O(t^{-\frac{1}{2}+\mbox{Im}(\vartheta(\lambda_{1}))})\\
    O(t^{-\frac{1}{2}-\mbox{Im}(\vartheta(\lambda_{1}))})  &  O(t^{-\frac{1}{2}})\\
\end{array}\right).
\end{align}

\section*{Appendix B}
According to Eq.\eqref{40}, we obtain
\begin{align}\label{80}
J^{(err)}=(M_{-}^{(err)})^{-1}M_{+}^{(err)}=M_{-}^{(app)}J^{(2)}(J^{(app)})^{-1}(M_{-}^{(app)})^{-1}.
\end{align}
Since $M_{-}^{(app)}=I$ on $\Sigma^{(err)}$, one has
\begin{align}\label{81}
J^{(err)}=J^{(2)}(J^{(app)})^{-1}.
\end{align}
Furthermore, we get
\begin{align}\label{82}
&J_{i}^{(err)}=J_{i}^{(2)}\ (i=1,2,\cdots, 8),\quad \mbox{outside} \quad D_{\lambda_{0}}^{\epsilon}\cup D_{\lambda_{1}}^{\epsilon},\notag\\
&J_{\lambda_{0}}^{(err)}=(J_{\lambda_{0}}^{(app)})^{-1}\qquad \mbox{on} \quad D_{\lambda_{0}}^{\epsilon},\notag\\
&J_{\lambda_{1}}^{(err)}=(J_{\lambda_{1}}^{(app)})^{-1}\qquad  \mbox{on} \quad D_{\lambda_{1}}^{\epsilon}.
\end{align}
Now, we will estimate the error of $J_{7}^{(err)}$ outside $D_{\lambda_{1}}^{\epsilon}$. On the jump contour $\lambda_{1}+\lambda_{1}\rho e^{\frac{3\pi i}{4}}(\rho>\epsilon)$, the jump matrix $J_{7}^{(err)}$ is
\begin{align}\label{83}
J_{7}^{(err)}=J_{7}^{(2)}=\left(\begin{array}{cc}
    1 &  r_{2}(\lambda)e^{-2ift}\delta_{+}^{2}\\
    0  &  1\\
\end{array}\right).
\end{align}
Observing
\begin{align}\label{84}
|e^{-2ift}|=e^{4\lambda_{1}^{2}\rho^{2}(\sqrt{2}\beta\lambda_{1}\rho-6\beta\lambda_{1}-\alpha)t}\leq e^{-\tilde{C}t},\ \tilde{C}>0.
\end{align}
Thus, we have
\begin{align}\label{85}
J_{7}^{(err)}-I=J_{7}^{(2)}-I=O(e^{-\tilde{C}t}).
\end{align}
Similarly, we have
\begin{align}\label{86}
J_{i}^{(err)}-I=J_{i}^{(2)}-I=O(e^{-\tilde{C}t}),\ i=1, 2\cdots, 8.
\end{align}
Moreover, from Eq.\eqref{78} and Eq.\eqref{79},  it is easy to find
\begin{align}\label{87}
J_{\lambda_{0}}^{(err)}-I=(J_{\lambda_{0}}^{(app)})^{-1}-I=\left(\begin{array}{cc}
    O(t^{-\frac{1}{2}}) &  O(t^{-\frac{1}{2}-\mbox{Im}(\vartheta(\lambda_{0}))})\\
    O(t^{-\frac{1}{2}+\mbox{Im}(\vartheta(\lambda_{0}))})  &  O(t^{-\frac{1}{2}})\\
\end{array}\right),\notag\\
J_{\lambda_{1}}^{(err)}-I=(J_{\lambda_{1}}^{(app)})^{-1}-I=\left(\begin{array}{cc}
    O(t^{-\frac{1}{2}}) &  O(t^{-\frac{1}{2}+\mbox{Im}(\vartheta(\lambda_{1}))})\\
    O(t^{-\frac{1}{2}-\mbox{Im}(\vartheta(\lambda_{1}))})  &  O(t^{-\frac{1}{2}})\\
\end{array}\right).
\end{align}

\section*{Appendix C}
The Cauchy integral formula on  contour $\Sigma$ can be defined as
\begin{align}\label{88}
(C_{\Sigma}(f))(\lambda)=\frac{1}{2\pi i}\int_{\Sigma}\frac{f(s)}{s-\lambda}ds.
\end{align}
Let
\begin{align}\label{89}
C^{-}_{V}(f)=C^{-}_{\Sigma}(f(V-I)),
\end{align}
where $V$ is a matrix given in $\Sigma$, and
$C^{+}_{\Sigma}, C^{+}_{\Sigma}$ denote the nontangential limits of the bounded operator $C_{\Sigma}$ approaching $\Sigma$ from left and right, respectively.

According to the RH problem  \eqref{41}, we can obtain
\begin{align}\label{90}
M^{(err)}-I=C_{\Sigma^{(err)}}M_{-}^{(err)}(J^{(err)}-I)=-\frac{1}{2\pi i\lambda}\int_{\Sigma^{(err)}}M_{-}^{(err)}(J^{(err)}-I)ds+O(\lambda^{-2}),
\end{align}
which indicates
\begin{align}\label{91}
M_{1}^{(err)}=-\frac{1}{2\pi i}\int_{\Sigma^{(err)}}M_{-}^{(err)}(J^{(err)}-I)ds.
\end{align}
Using Holder inequality, we have
\begin{align}\label{92}
\mid M_{1}^{(err)}\mid\leq C_{1}||M_{-}^{(err)}-I||_{L^{2}}||J^{(err)}-I||_{L^{2}}+C_{2}||J^{(err)}-I||_{L^{1}}, \ C_{1}, C_{2}>0.
\end{align}
Beside, it is not hard to get $||M_{-}^{(err)}-I||_{L^{2}}\leq C_{3}||J^{(err)}-I||_{L^{2}}, C_{3}>0$, then we finally obtain
\begin{align}\label{93}
\mid M_{1}^{(err)}\mid\leq C_{1}C_{3}||J^{(err)}-I||_{L^{2}}+C_{2}||J^{(err)}-I||_{L^{1}}.
\end{align}
Combining Eq.\eqref{86} and Eq.\eqref{87}, we arrive at
\begin{align}\label{94}
|M_{1}^{(err)}(x, t, \lambda)|=
O(t^{-\frac{1}{2}-\mbox{max}\{|\mbox{Im}\vartheta(\lambda_{0})|,|\mbox{Im}\vartheta(\lambda_{1})|\}}).
\end{align}

\section*{Appendix D}
The solution $M^{(mod)}(\tilde{\lambda})$ of the model RH problem \eqref{54} can be given explicitly via using the  Liouville's theorem and parabolic cylinder functions. Since the jump matrix $J^{(mod)}$  is constant, the logarithmic derivative $\frac{d}{d\tilde{\lambda}}M^{(mod)}(M^{(mod)})^{-1}$ possesses continuous jump along any of the rays, which indicates that $M^{(mod)}$  solves the following ordinary differential equation
\begin{align}\label{69}
\frac{d}{d\tilde{\lambda}}M^{(mod)}+\left(\begin{array}{cc}
    -\frac{i}{2}\tilde{\lambda}  &  \Psi\\
  \Phi &  \frac{i}{2}\tilde{\lambda}\\
\end{array}\right)M^{(mod)}=0,
\end{align}
where $\Psi=i\left[M_{1}^{(\infty)}\right]_{12}, \Phi=-i\left[M_{1}^{(\infty)}\right]_{21}$.
The solution of $\eqref{69}$ can be written as
\begin{align}\label{70}
M^{(mod)}=\left(\begin{array}{cc}
    M_{11}^{(mod)}  &  \frac{-\frac{i}{2}\tilde{\lambda}M_{22}^{(mod)}-\frac{dM_{22}^{(mod)}}{d\tilde{\lambda}}}{\Phi}\\
  \frac{\frac{i}{2}\tilde{\lambda}M_{11}^{(mod)}-\frac{dM_{11}^{(mod)}}{d\tilde{\lambda}}}{\Psi} &  M_{22}^{(mod)}\\
\end{array}\right),
\end{align}
where the functions $M_{jj}^{(mod)}, j=1, 2,$ satisfy the parabolic cylinder equations
\begin{align}\label{71}
\frac{d^{2}}{d\tilde{\lambda}^{2}}M_{11}^{(mod)}+(-\frac{i}{2}-\Phi\Psi+\frac{\tilde{\lambda}^{2}}{4})M_{11}^{(mod)}=0,\notag\\
\frac{d^{2}}{d\tilde{\lambda}^{2}}M_{22}^{(mod)}-(-\frac{i}{2}+\Phi\Psi-\frac{\tilde{\lambda}^{2}}{4})M_{22}^{(mod)}=0.
\end{align}
According to the property of standard parabolic cylinder equation and $M_{11}^{(mod)}\rightarrow e^{\frac{1}{4}i\tilde{\lambda}^{2}}\tilde{\lambda}^{-i\vartheta},$ $ M_{22}^{(mod)}\rightarrow e^{-\frac{1}{4}i\tilde{\lambda}^{2}}\tilde{\lambda}^{i\vartheta}, \tilde{\lambda}\rightarrow\infty$, we obtain
\begin{align}\label{72}
M_{11}^{(mod)}=\left\{
\begin{array}{lr}
(ie^{-\frac{3\pi}{4}i})^{i\vartheta}D_{-i\vartheta}(ie^{-\frac{3\pi}{4}i}\tilde{\lambda}) \qquad \mbox{Im}(\tilde{\lambda})>0,\\
(ie^{\frac{\pi}{4}i})^{i\vartheta}D_{-i\vartheta}(ie^{\frac{\pi}{4}i}\tilde{\lambda}) \qquad \qquad \mbox{Im}(\tilde{\lambda})<0,
  \end{array}
\right.
\end{align}
\begin{align}\label{73}
M_{22}^{(mod)}=\left\{
\begin{array}{lr}
(ie^{-\frac{\pi}{4}i})^{-i\vartheta}D_{i\vartheta}(ie^{-\frac{\pi}{4}i}\tilde{\lambda}) \ \qquad \mbox{Im}(\tilde{\lambda})>0,\\
(ie^{\frac{3\pi}{4}i})^{-i\vartheta}D_{i\vartheta}(ie^{\frac{3\pi}{4}i}\tilde{\lambda}) \quad \qquad \mbox{Im}(\tilde{\lambda})<0.
  \end{array}
\right.
\end{align}
Then, we can get
\begin{gather}
M^{(mod)}_{-}(\tilde{\lambda})^{-1}M^{(mod)}_{+}(\tilde{\lambda})=M^{(mod)}_{-}(0)^{-1}M^{(mod)}_{+}(0)=\notag\\
\left(\begin{array}{cc}
  (ie^{\frac{\pi}{4}i})^{i\vartheta}\frac{2^{\frac{-i\vartheta}{2}}\sqrt{\pi}}{\Gamma(\frac{1+i\vartheta}{2})}    & (ie^{\frac{3\pi}{4}i})^{1-i\vartheta}\frac{2^{\frac{1+i\vartheta}{2}}\sqrt{\pi}}{\Phi\Gamma(\frac{-i\vartheta}{2})} \\
  (ie^{\frac{\pi}{4}i})^{1+i\vartheta}\frac{2^{\frac{1-i\vartheta}{2}}\sqrt{\pi}}{\Psi\Gamma(\frac{i\vartheta}{2})} &   (ie^{\frac{3\pi}{4}i})^{-i\vartheta}\frac{2^{\frac{i\vartheta}{2}}\sqrt{\pi}}{\Gamma(\frac{1-i\vartheta}{2})}\\
\end{array}\right)^{-1}\notag\\
\left(\begin{array}{cc}
  (ie^{-\frac{3\pi}{4}i})^{i\vartheta}\frac{2^{\frac{-i\vartheta}{2}}\sqrt{\pi}}{\Gamma(\frac{1+i\vartheta}{2})}    & (ie^{-\frac{\pi}{4}i})^{1-i\vartheta}\frac{2^{\frac{1+i\vartheta}{2}}\sqrt{\pi}}{\Phi\Gamma(\frac{-i\vartheta}{2})} \\
  (ie^{-\frac{3\pi}{4}i})^{1+i\vartheta}\frac{2^{\frac{1-i\vartheta}{2}}\sqrt{\pi}}{\Psi\Gamma(\frac{i\vartheta}{2})} &   (ie^{-\frac{\pi}{4}i})^{-i\vartheta}\frac{2^{\frac{i\vartheta}{2}}\sqrt{\pi}}{\Gamma(\frac{1-i\vartheta}{2})}\\
\end{array}\right)\notag\\
=\left(\begin{array}{cc}
    1-r_{1}(\lambda_{1})r_{2}(\lambda_{1})  &  -r_{2}(\lambda_{1})\\
  r_{1}(\lambda_{1}) &  1\\
\end{array}\right),
\end{gather}
which leads to
\begin{align}\label{74}
\Psi=\frac{\sqrt{2\pi}(-1)^{i\vartheta(\lambda_{0})+\frac{1}{2}}e^{-\frac{3\pi i}{4}+\frac{\pi\vartheta(\lambda_{0}) }{2}}}{r_{1}(\lambda_{0})\Gamma(i\vartheta(\lambda_{0}))},\quad \Phi=-\frac{\sqrt{2\pi}(-1)^{-i\vartheta(\lambda_{0})+\frac{1}{2}}e^{-\frac{\pi i}{4}+\frac{\pi\vartheta(\lambda_{0}) }{2}}}{r_{2}(\lambda_{0})\Gamma(-i\vartheta(\lambda_{0}))}.
\end{align}

\section*{Acknowledgements}
\hspace{0.3cm}
This work was supported by the project is supported by National Natural Science Foundation of China(No.12175069)
and Science and Technology Commission of Shanghai Municipality (No.21JC1402500 and
No.18dz2271000).

\end{document}